\documentclass{article}%
\usepackage{amsmath}
\usepackage{amsfonts}
\usepackage{amssymb}
\usepackage{graphicx}%
\providecommand{\U}[1]{\protect\rule{.1in}{.1in}}

\begin{document}

\title{The Concept of Statistical Evidence\\{\large Historical Roots and Current Developments}}
\author{Michael Evans\\Department of Statistical Sciences, University of Toronto}
\date{}
\maketitle

\begin{abstract}
One can argue that one of the main roles of the subject of statistics is to
characterize what the evidence in collected data says about questions of
scientific interest. There are two broad questions that we will refer to as
the estimation question and the hypothesis assessment question. For
estimation, the evidence in the data should determine a particular value of an
object of interest together with a measure of the accuracy of the estimate,
while for hypothesis assessment, the evidence in the data should provide
evidence in favor of or against some hypothesized value of the object of
interest together with a measure of the strength of the evidence. This will be
referred to as the evidential approach to statistical reasoning which can be
contrasted with the behavioristic or decision-theoretic approach where the
notion of loss is introduced and the goal is to minimize expected losses.
While the two approaches often lead to similar outcomes, this is not always
the case and it is commonly argued that the evidential approach is more suited
to scientific applications. This paper traces the history of the evidential
approach and summarizes current developments.

\end{abstract}

\section{Introduction}

Most statistical analyses refer to the concept of \textit{statistical
evidence} as in phrases like "the evidence in the data suggests" or "based on
the evidence we conclude", etc. It has long been recognized, however, that the
concept itself has never been satisfactorily defined or, at least no
definition has been offered that has met with general approval. This article
is about the concept of statistical evidence outlining some of the key
historical aspects of the discussion and the current state of affairs.

One issue that needs to be dealt with right away is whether or not it is even
necessary to settle on a clear definition. After all Statistics has been
functioning as an intellectual discipline for many years without a resolution.
There are at least two reasons why resolving this is important.

First, to be vague about the concept leaves open the possibility of
misinterpretation and ambiguity as in "if we don't know what statistical
evidence is, how can we make any claims about what the evidence is saying?"
This leads to a degree of ad hocracy in the subject where different analysts
measure and interpret the concept in different ways. For example, consider the
replicability of research findings in many fields where statistical
methodology is employed. One can't claim that being more precise about
statistical evidence will fix such problems, but it is reasonable to suppose
that establishing a sound system of statistical reasoning based on a clear
prescription of what we mean by this concept, can only help.

Second, the subject of Statistics cannot claim that it speaks with one voice
on what constitutes a correct statistical analysis. For example, there is the
Bayesian versus frequentist divide as well as the split between the evidential
versus the decision-theoretic or behavioristic approaches to determining
inferences. This diversity of opinion, while interesting in and of itself,
does not enhance general confidence in the soundness of the statistical
reasoning process.

As such, it is reasonable to argue that settling the issue of what statistical
evidence is and how it is to be used to establish inferences is of paramount
importance. Statistical reasoning is a substantial aspect of how many
disciplines, from Anthropology to Particle and Quantum Physics and on to
Zoology, determine truth in their respective subjects. One can claim that the
role of the subject of Statistics is to provide these scientists with a system
of reasoning that is logical, consistent and sound in the sense that it
produces satisfactory results in practical contexts free of paradoxes. Recent
controversies over the use of p-values, that have arisen in a number of
scientific contexts, suggest that, at the very least, this need is not being met.

In Section 2 we provide some background. Section 3 discusses various attempts
at measuring statistical evidence and why these are not satisfactory. Section
4 discusses what we call the Principle of Evidence and how it can resolve many
of the difficulties associated with defining and using the concept of
statistical evidence.

\section{The Background}

There are basically two problems that the subject of Statistics has been
developed to address which we will refer to as the estimation (\textbf{E})
problem and the hypothesis assessment (\textbf{H}) problem. We will denote the
data produced by some process by the symbol $x$ and always assume that this
has been collected correctly, namely, by random sampling from relevant
populations. Obviously, this requirement is not always met in applications but
the focus here is on the ideal situation that can serve as a guide in
practical contexts. The fact that data have not been collected properly in an
application can be, however, a substantial caveat when interpreting the
results of a statistical analysis and at least should be mentioned as a qualification.

Now suppose there is an object or concept of interest, denoted by $\Psi,$ in
the world whose value is unknown. For example, $\Psi$ could be the half-life
of some fundamental particle, the rate of recidivism in a certain jurisdiction
at some point in time, the average income and assets of voters in a given
country, a graph representing the relationships among particular variables
measured on each element of a population, etc. A scientist collects data $x$
which they believe contains evidence concerning the answers to the following questions.

\begin{quote}
\textbf{E}: what value does $\Psi$ take (estimation)?\smallskip\ 

\textbf{H}: does $\Psi=\psi_{0}$ (hypothesis assessment)?
\end{quote}

\noindent Actually, there is more to these questions than just providing
answers. For example, when quoting an estimate of $\Psi$ based on $x,$ there
should also be an assessment of how accurate it is believed that the estimate
is. Also, if a conclusion is reached concerning whether or not $\Psi=\psi_{0}$
based on $x,$ there should also be an assessment of how strong a basis there
is for drawing this conclusion. At least as far as this article goes,
Statistics as a subject is concerned with \textbf{E} and \textbf{H}. For
example, if interest is in predicting a future value of some process based on
observing $x,$ then $\Psi$ will correspond to this future value.

The difference between the evidential and decision-theoretic approach can now
be made clear. The evidential approach is concerned with summarizing what the
evidence in $x$ says about the answers to \textbf{E} and \textbf{H}. For
example, for \textbf{E} a conclusion might be that the evidence in the data
implies that the true value of $\Psi$ is best estimated by $\psi(x)$ and also
provide an assessment of the accuracy of $\psi(x)$. For \textbf{H} the
conclusion would be that there is evidence against/in favor of $\Psi=\psi_{0}$
together with a measure of the strength of this evidence. By contrast the
decision-theoretic approach adds the idea of utility or loss to the problem
and searches for answers which optimize in a statistical sense expected
utility/loss as in quoting an optimal estimate $\psi(x)$ or optimally
accepting/rejecting $\Psi=\psi_{0}.$ In a sense the decision-theoretic
approach dispenses with the idea of evidence and tries to characterize optimal
behavior. A typical complaint about the decision-theoretic approach is that
for scientific problems, where the goal is the truth whatever that might be,
there is no role for the concept of utility/loss. That is not to say that
there is no role for decision-theory in applications generally, but for
scientific problems we want to know as clearly as possible what the evidence
in the data is saying. There may indeed be good reasons to do something that
contradicts the evidence, now being dictated by utility considerations, but
there is no reason not to state what the evidence says and then justify the
difference. By and large we only consider the evidential approach here but it
will be seen that there is some confounding at times between the two and this
needs to be discussed.

The other significant division between the approaches to \textbf{E} and
\textbf{H} is between frequentist and Bayesian statistics. In fact, these are
also subdivisions of the evidential and the decision-theoretic approach. As
they do not generally lead to the same results, this division needs to be
addressed as part of our discussion. As a broad characterization, frequentism
considers a thought experiment where it is imagined that the observed data $x$
could be observed as a member of a sequence of basically identical situations,
say obtaining $x_{1},x_{2},\ldots,$ and then how well an inference does in
such sequences, in terms of reflecting the truth concerning $\Psi,$ is
measured. Inferences with good performance are considered more reliable and
thus preferred. As soon as these performance measures are invoked, however,
this brings things close to being part of decision theory and, as will be
discussed here, there doesn't seem to be a clear way to avoid this. By
contrast the Bayesian approach insists that all inferences must be dictated by
the observed data without taking into account hypothetical replications as a
measure of reliability. As will be discussed in Section 3.4.1, however, it is
possible to resolve this apparent disagreement to a great extent.

\section{Evidential Inference}

Statistical evidence is the central core concept of this approach. It is now
examined how this has been addressed by various means and to what extent these
have been successful. For this we need some notation. The basic ingredient of
every statistical theory is the sampling model $\{f_{\theta}:\theta\in
\Theta\},$ a collection of probability density functions on some sample space
$\mathcal{X}$ with respect to some support measure $\nu.$ The variable
$\theta$ is the \textit{model parameter} which takes values in the
\textit{parameter space} $\Theta$ and it indexes the probability distributions
in the model. The idea is that one of these distributions in the model
produced the data $x$ and this is represented as $f_{\theta_{true}}.$ The
object of interest $\Psi$ is represented as a function $\Psi:\Theta
\rightarrow\Psi(\Theta)$ where $\psi=\Psi(\theta)$ can be obtained from the
distribution $f_{\theta}.$ The questions \textbf{E} and \textbf{H} can
be\textbf{ }answered categorically if $\psi_{true}=\Psi(\theta_{true})$
becomes known.

A simple example illustrates these concepts and will be used throughout this
discussion.\smallskip

\noindent\textbf{Example} \textbf{1.} \textit{location normal.}

Suppose $x=(x_{1},\ldots,x_{n})$ is an $iid$ sample of $n$ from a distribution
in the family $\{n(\cdot,\mu,\sigma_{0}^{2}):\mu\in%
\mathbb{R}
^{1}\}$, where $n(\cdot,\mu,\sigma_{0}^{2})$ denotes the density of a normal
distribution with unknown mean $\mu$ and known variance $\sigma_{0}^{2},$ and
the goal is inference about the true value of $\psi=\Psi(\mu)=\mu.$ While
$\Psi$ here is one-to-one it could be many-to-one, as in $\psi=\Psi(\mu
)=|\mu|,$ and such contexts raise the important issue of nuisance parameters
and how to deal with them$.$ $\blacksquare$

\subsection{P-values, E-values and Confidence Region}

The \textit{p-value} is the most common statistical approach to measuring
evidence associated with \textbf{H}. The p-value is associated with Fisher
(1925) although there are precursors. The history of this concept is discussed
in Stigler (1986).

We suppose that there is a statistic $T_{\psi_{0}}(x)$ and consider its
distribution under the hypothesis $H_{0}:\Psi(\theta)=\psi_{0}.$ The idea
behind the p-value is to answer the question: is the observed value of
$T_{\psi_{0}}(x)$ something we would expect to see if $\psi_{0}$ is the true
value of $\Psi?$ If $T_{\psi_{0}}(x)$ is a surprising value, then this is
interpreted as evidence against $H_{0}.$ A p-value measures the location of
$T_{\psi_{0}}(x)$ in the distributions of $T_{\psi_{0}}$ under $H_{0}$ with
small values of the p-value being an indication that the observed value is
surprising. For example, it is quite common that large values of $T_{\psi_{0}%
}(x)$ are in the tails (regions of low probability content) of each of the
distributions of $T_{\psi_{0}}$ under $H_{0}$ and so, given $T_{\psi_{0}}(x),$%
\[
p_{\psi_{0}}(x)=\sup_{\Psi(\theta)=\psi_{0}}P_{\theta}(T_{\psi_{0}}(X)\geq
T_{\psi_{0}}(x))
\]
is computed as the \textit{p-value}. \smallskip

\noindent\textbf{Example} \textbf{2.} \textit{location normal.}

Suppose that it is required to assess the hypothesis $H_{0}:\mu=\mu_{0}.$ For
this it is common to use the statistic $T_{\mu_{0}}(x)=\sqrt{n}|\bar{x}%
-\mu_{0}|/\sigma_{0},$ as this has a fixed distribution (the absolute value of
a standard normal variable) under $H_{0},$ with the p-value given by
\[
p_{_{\mu_{0}}}(x)=P_{_{\mu_{0}}}(T_{\mu_{0}}(X)\geq T_{\mu_{0}}(x))=2(1-\Phi
(\sqrt{n}|\bar{x}-\mu_{0}|/\sigma_{0}))
\]
where $\Phi$ is the $N(0,1)$ cdf. $\blacksquare$\smallskip

Several issues arise with the use of p-values in general. First, what value
$\alpha$ is small enough to warrant the assertion that evidence against
$H_{0}$ has been obtained when the p-value is less than $\alpha?$ It is quite
common in many fields to use the value $\alpha=0.05$ as a cut-off but this is
not universal as in particle physics $\alpha=5.733031\times10^{-7}$ is
commonly employed. Recently, due to concerns with replication of results, it
has been suggested that the value $\alpha=0.005$ be used as a standard. The
issue of the appropriate cut-off is not resolved.

It is also not the case that a value greater than $\alpha$ is to be
interpreted as evidence in favor of $H_{0}$ being true. For note that, in the
case that $T_{\psi_{0}}$ has a single continuous distribution $P_{H_{0}}$
under $H_{0},$ then $P_{H_{0}}(T_{\psi_{0}}\geq T_{\psi_{0}}(x))$ has a
uniform$(0,1)$ when $x\sim P_{H_{0}}.$ This implies that a p-value near 0 has
the same probability of occurring as a p-value near 1. It is often claimed
that a valid p-value must have this uniformity property. But consider the
p-value of Example 3 where $\bar{x}\sim N(\mu_{0},\sigma_{0}^{2}/n)$ under
$H_{0},$ and so, as $n$ rises, the distribution of $\bar{x}$ becomes more and
more concentrated about $\mu_{0}.$ For large enough $n$ virtually all of the
distribution of $\bar{x}$ under $H_{0}$ will be concentrated in the interval
$(\mu_{0}-\delta,\mu_{0}+\delta)$ where $\delta$ represents a deviation from
$\mu_{0}$ that is of scientific interest while a smaller deviation is not,
e.g., the measurements $x_{i}$ are made to an accuracy no greater than
$2\delta.$ Note that in any scientific problem it is necessary to state the
accuracy with which the investigator wants to know the object $\Psi$ as this
guides the measurement process as well as the choice of sample size. The
p-value ignores this fact and, when $n$ is large enough, could record evidence
against $H_{0}$ when in fact the data is telling us that $H_{0}$ is
effectively true and evidence in favor should be stated. This distinction
between statistical significance and scientific significance has long been
recognized, see Boring (1919), and needs to be part of statistical
methodology, see Example 5. This fact also underlies a recommendation,
although not commonly followed as it doesn't really address the problem, that
the $\alpha$ cut-off should be reduced as the sample size increases. Perhaps
the most important take-away from this is that a value $T_{\psi_{0}}(x)$ that
lies in the tails of its distributions under $H_{0}$ is not necessarily
evidence against $H_{0}.$ It is commonly stated that a confidence region for
$\Psi$ should also be provided but this does not tell us anything about
$\delta,$ which needs to be given as part of the problem.

There are many other issues that can be raised about the p-value where many of
these are covered by the word \textit{p-hacking} and associated with the
choice of $\alpha.$ For example, as discussed in Cornfield (1966), suppose an
investigator is using a particular $\alpha$ as the cut-off for evidence
against and, based on a data set of size $n_{1}$ obtains the p-value
$\alpha+\epsilon$ where $\epsilon$ is small. Since finding evidence against
$H_{0}$ is often regarded as a positive, as it entails a new discovery, the
investigator decides to collect an additional $n_{2}$ data values, computes
the p-value based on the full $n=n_{1}+n_{2}\ $data values and obtains a
p-value less than $\alpha.$ But this ignores the two stage structure of the
investigation and, when this is taken into account, the probability of finding
evidence against $H_{0}$ at\ level $\alpha$ when $H_{0}$ is true, and assuming
a single distribution under $H_{0},$ equals $P_{H_{0}}(A)+P_{H_{0}}(A^{c}\cap
B)$ where $A=$ $\{$evidence against $H_{0}$ found at first stage$\}$ and $B=$
$\{$evidence against $H_{0}$ found at second stage$\}$. If $P_{H_{0}%
}(A)=\alpha,$ then $P_{H_{0}}(A)+P_{H_{0}}(A^{c}\cap B)>\alpha.$ So, even
though the investigator has done something very natural, the logic of the
statistical procedure based on using p-values with a cut-off, effectively
prevents finding evidence against $H_{0}.$

Royall (1997) provides an excellent discussion of the deficiencies of p-values
and \textit{rejection trials (}p-values with cut-off $\alpha)$ when
considering these as measuring evidence\textit{. }Sometimes it is recommended
that the p-value itself be reported without reference to a cut-off $\alpha$
but we have seen already that a small p-value does not necessarily constitute
evidence against and so this is not a solution.

Confidence regions are intimately connected with p-values. Suppose, for each
$\psi_{0}\in\Psi(\Theta)$ there is a statistic $T_{\psi_{0}}$ that produces a
valid p-value for $H_{0}$ as has been described and an $\alpha$ cut-off is
used. Now put $C_{1-\alpha}(x)=\{\psi_{0}:p_{\psi_{0}}(x)>\alpha\}.$ Then
$P_{\theta}(\Psi(\theta)\in C_{1-\alpha}(x))\geq1-\alpha$ and so $C_{1-\alpha
}(x)$ is a $(1-\alpha)$\textit{-confidence region} for $\psi.$ As such,
$\psi_{0}\in C_{1-\alpha}(x)$ is equivalent to $\psi_{0}$ not being rejected
at level $\alpha$ and so all the problematical issues with p-values as
measures of evidence apply equally to confidence regions. Moreover, it is
correctly stated that, when $x$ has been observed, then $1-\alpha$ is not a
lower bound on the probability that $\psi_{true}\in C_{1-\alpha}(x).$ Of
course that is what we want, namely, to state a probability that measures our
belief that $\psi_{true}$ is in this set and so confidence regions are
commonly misinterpreted, however, see the discussion of bias in Section 3.4.1.

An alternative approach to using p-values is provided by e-values, see Shafer
et al. (2011), Vovk and Wang (2023) and Grunwald et al. (2024). An
\textit{e-variable} for a hypothesis $H_{0}:\Psi(\theta)=\psi_{0}$\ is a
nonnegative statistic $N_{\psi_{0}}$ that satisfies $E_{\theta}(N_{\psi_{0}%
})\leq1$ whenever $\Psi(\theta)=\psi_{0}.$ The observed value $N_{\psi_{0}%
}(x)$ is called an \textit{e-value }where the "e" stands for expectation to
contrast it with the "p" in p-value which stands for probability. Also, a
cut-off $\alpha\in(0,1)$ needs to be specified such that $H_{0}$ is rejected
whenever $N_{\psi_{0}}(x)\geq1/\alpha.$ It is immediate, from Markov's
inequality, that $P_{\theta}(N_{\psi_{0}}(x)\geq1/\alpha)\leq\alpha$ whenever
$\Psi(\theta)=\psi_{0}$ and so this provides a rejection trial with cut-off
$\alpha$ for $H_{0}.$\smallskip

\noindent\textbf{Example} \textbf{3.} \textit{location normal.}

For the context of Example 2 define $N_{\mu_{0}}(x)=a[2(1-\Phi(\sqrt{n}%
|\bar{x}-\mu_{0}|/\sigma_{0}))]^{a-1}$ for any $a\in(0,1)$ and it is immediate
that $N_{\mu_{0}}$ is an e-variable for $H_{0}:\mu=\mu_{0}.$ $\blacksquare
$\smallskip

Consider the situation where data is collected sequentially, $N_{i,\psi_{0}%
}(x_{i})$ is an e-variable for $H_{0}:\Psi(\theta)=\psi_{0}$ based on
independent data $x_{i}$ and there is a stopping time $S,$ e.g., stop when
$N_{\psi_{0}}(x_{1},\ldots,x_{n})=%
{\displaystyle\prod_{i=1}^{n}}
N_{i,\psi_{0}}(x_{i})>1/\alpha.$ Also, put $N_{0,\psi_{0}}\equiv1.$ Then,
whenever $\Psi(\theta)=\psi_{0},$ we have that%
\begin{align*}
E_{\theta}(N_{\psi_{0}}(x_{1},\ldots,x_{n})\,|\,x_{1},\ldots,x_{n-1})  &
=N_{\psi_{0}}(x_{1},\ldots,x_{n-1})E_{\theta}(N_{n,\psi_{0}}(x_{n}))\\
&  \leq N_{\psi_{0}}(x_{1},\ldots,x_{n-1})
\end{align*}
and so the process $N_{\psi_{0}}(x_{1},\ldots,x_{n})$ is a discrete time
super-martingale. Moreover, it is clear that $N_{\psi_{0}}(x_{1},\ldots
,x_{n})$ is an e-variable for $H_{0}:\Psi(\theta)=\psi_{0}.$ This implies
that, under conditions, $E_{\theta}(N_{\psi_{0}}(x_{1},\ldots,x_{S}))\leq1$
and so the stopped variable $N_{\psi_{0}}(x_{1},\ldots,x_{S})$ is also an
e-variable for $H_{0}:\Psi(\theta)=\psi_{0}$. Assuming stopping time
$S_{\alpha}=\inf\{n:N_{\psi_{0}}(x_{1},\ldots,x_{n})\geq1/\alpha\}$ is finite
with probability 1 when $H_{0}$ holds, then by Ville's inequality $P\left(
N_{\psi_{0}}(x_{1},\ldots,x_{S_{\alpha}})\geq1/\alpha\right)  \leq\alpha.$
This implies that the problem for p-values when sampling until rejecting
$H_{0}$ at size $\alpha,$ is avoided when using e-values.

While e-values have many interesting and useful properties relative to
p-values, the relevant question here is whether or not these serve as measures
of statistical evidence. Usage of both requires the specification of $\alpha$
to determine when there are grounds for rejecting $H_{0}.$ Sometimes this is
phrased instead as "evidence against $H_{0}$ has been found" but, given the
arbitrariness of the choice of $\alpha$ and the failure to properly express
when evidence in favor of $H_{0}$ has been found, neither seems suitable as an
expression of statistical evidence. One\ could argue that the intention behind
these approaches is not to characterize statistical evidence but rather to
solve the reject/accept problems of decision theory. It is the case, however,
at least for p-values, that these are used as if they are proper
characterizations of evidence and this does not seem suitable for purely
scientific applications.

Another issue that needs to be addressed in a problem is how to find the
statistic $T_{\psi_{0}}$ or $N_{\psi_{0}}.$ A common argument is to use
likelihood ratios, see Sections 3.2 and 3.3, but this does not resolve the
problems that have been raised here and there are difficulties with the
concept of likelihood that need to be addressed as well.

\subsection{Birnbaum on Statistical Evidence}

Alan Birnbaum devoted considerable thought to the concept of statistical
evidence in the frequentist context. Birnbaum (1962) contained what seemed
like a startling result about the implications to be drawn from the concept
and his final paper Birnbaum (1977) contained a proposal for a definition of
statistical evidence. Also, see Giere (1977) for a full list of Birnbaum's
publications many of which contain considerable discussion concerning
statistical evidence.

Birnbaum (1962) considered various relations, as defined by statistical
principles, on the set $\mathcal{I}$ of all inference bases where an
\textit{inference base} is the pair $I=(\{f_{\theta}:\theta\in\Theta\},x)$
consisting of a sampling model and data supposedly generated from a
distribution in the model. A \textit{statistical principle} is then a relation
defined on $\mathcal{I},$ namely, a subset of $\mathcal{I\times I}.$ There are
three basic statistical principles that are commonly invoked as part of
evidential inference, namely, the likelihood principle $(L),$ the sufficiency
principle $(S)$ and the conditionality principle $(C).$ Inference bases
$I_{i}=(\{f_{i\theta}:\theta\in\Theta\},x_{i}),$ for $i=1,2,$ satisfy
$(I_{1},I_{2})\in L$ if for some constant $c>0$ we have $f_{1\theta}%
(x_{1})=cf_{2\theta}(x_{2})$ for every $\theta\in\Theta,$ they satisfy
$(I_{1},I_{2})\in S$ if the models have equivalent (1-1 functions of each
other) \textit{minimal sufficient statistics} $m_{i}$ ( a sufficient statistic
that is a function of every other sufficient statistic) that take the
equivalent values at the corresponding data values and $(I_{1},I_{2})\in C$ if
there is \textit{ancillary statistic} $a$ (a statistic whose distribution is
independent of $\theta)$ such that $I_{2}$ can be obtained from $I_{1}$ (or
conversely) via conditioning on $a(x)$ so $I_{2}=(\{f_{1\theta}(\cdot
\,|\,a(x)):\theta\in\Theta\},x_{1}).$ The basic idea is that, if two inference
bases are related by one of these principles, then they contain the same
statistical evidence concerning concerning the unknown true value of $\theta.$
For this idea to make sense, it must be the case that these principles form
equivalence relations on $\mathcal{I}.$ In Evans (2013) it is shown that $L$
and $S$ do form equivalence relations but $C$ does not and this latter result
is connected with the fact that a unique \textit{maximal ancillary} (an
ancillary which every other ancillary is a function of) does not exist.

Birnbaum (1962) provided a proof of the result, known as Birnbaum's Theorem,
that if a statistician accepted the principles $S$ and $C,$ then they must
accept $L.$ This is highly paradoxical because frequentist statisticians
generally accept both $S$ and $C$ but not $L$ as $L$ does not permit repeated
sampling. Two very different sampling models can lead to proportional
likelihood functions so repeated sampling behavior is irrelevant under $L$.
There has been considerable discussion over the years concerning the validity
of this proof but no specific flaw was found. A resolution of this is provided
in Evans (2013) where it is shown that $S\cup C$ is not an equivalence
relation and what Birnbaum actually proved is that the smallest equivalence
relation on $\mathcal{I}$ that contains $S\cup C$ is $L.$ This substantially
weakens the result as there is no reason to accept the additional generated
equivalences and in fact it makes more sense to consider the largest
equivalence relation on $\mathcal{I}$ that is contained in $S\cup C.$ Also, in
Evans, Fraser and Monette (1986), it was shown, via an argument similar to the
one found in Birnbaum (1962), that a statistician who accepts $C$ alone must
accept $L$ although, as shown in Evans (2013), this simply shows that the
smallest equivalence relation on $\mathcal{I}$ that contains $C$ is $L.$ As
shown in Evans and Frangakis (2023) issues concerning $C$ can be resolved by
restricting conditioning on ancillary statistics to \textit{stable
ancillaries} (those ancillaries such that conditioning on them retains the
ancillarity of all other ancillaries and similarly have their ancillarity
retained when conditioning on any other ancillary) as there always is a
maximal stable ancillary. This provides a conditionality principle that is a
proper characterization of statistical evidence but it does not lead to $L.$

While Birnbaum did not ultimately resolve the issues concerning statistical
evidence, Birnbaum (1977) makes a suggestion as a possible starting point by
proposing the \textit{confidence concept }in the context of comparing two
hypotheses $H_{0}$ versus $H_{1}.$ The confidence concept is characterized by
two \textit{error probabilities}
\begin{align*}
\alpha &  =\text{probability of rejecting (accepting) }H_{0}\text{ }%
(H_{1})\text{ when it is true (false),}\\
\beta &  =\text{probability of accepting (rejecting) }H_{0}\text{ }%
(H_{1})\text{ when it is false (true),}%
\end{align*}
and then reporting $(\alpha,\beta)$ with the following interpretation:%
\begin{align*}
&  \text{rejecting\ }H_{0}\text{ constitutes strong evidence against }%
H_{0}\text{ (in favor of }H_{1})\\
&  \text{when }\alpha\text{ and }\beta\text{ are small.}%
\end{align*}
This clearly results from a confounding of the decision theoretic (as in
Neyman-Pearson) approach to hypothesis testing with the evidential approach,
as rejection trials similarly do. Also, this does not give a general
definition of what is meant by statistical evidence as it really only applies
in the simple versus simple hypothesis testing context and it suffers from
many of the same issues as discussed concerning p-values. From reading
Birnbaum's papers it seems he had largely despaired of ever finding a fully
satisfactory characterization of statistical evidence in the frequentist
context. It will be shown in Section 3.4, however, that the notion of
frequentist error probabilities do play a key role in a development of the concept.

\subsection{Likelihood}

Likelihood is another statistical concept initiated by Fisher (1922). While
the likelihood function plays a key role in frequentism, Edwards (1992) and
Royall (1997) develop a theory of inference, called \textit{pure likelihood
theory}, based solely on the likelihood function. The basic axiom is thus the
likelihood principle $L$ which says that the likelihood function
$L(\theta\,|\,x)=cf_{\theta}(x),$ from inference base $I=(\{f_{\theta}%
:\theta\in\Theta\},x),$ completely summarizes the evidence given in $I$
concerning the true value of $\theta$ and repeated sampling does not play a
role$.$ In particular, the ratio $L(\theta_{1}\,|\,x)/$ $L(\theta_{2}\,|\,x)$
provides the evidence for $\theta_{1}$ relative to the evidence for
$\theta_{2}$ and this is independent of the arbitrary constant $c.$

While this argument for the ratios seems acceptable, a problem arises when we
ask: does the value $L(\theta\,|\,x)$ provide evidence in favor of or against
$\theta$ being the true value? To avoid the arbitrary constant Royall (1997)
replaces the likelihood function by the \textit{relative likelihood function
}given by\textit{ }$L^{rel}(\theta\,|\,x)=L(\theta\,|\,x)/L(\theta
_{MLE}(x)\,|\,x)$ where $\theta_{MLE}(x)=\arg\sup_{\theta}L(\theta\,|\,x)$ is
the \textit{maximum likelihood estimate}. The relative likelihood always takes
values in $[0,1].$ Note that, if $L^{rel}(\theta\,|\,x)\geq1/r,$ then
$L^{rel}(\theta^{\prime}\,|\,x)/L^{rel}(\theta\,|\,x)\leq r$ for all
$\theta^{\prime}$ and so no other value is supported by more than $r$ times
the support accorded to $\theta.$ For $H_{0}:\theta=\theta_{0},$ Royall (1997)
then argues, based on single urn experiment, for $L^{rel}(\theta
_{0}\,|\,x)\geq1/8$ to represent very strong evidence in favor and for
$L^{rel}(\theta_{0}\,|\,x)\geq1/32$ to represent quite strong evidence in
favor of $\theta.$ Certainly, these values seem quite arbitrary and again, as
with p-values, we don't have a clear cut-off between evidence in favor and
evidence against. Whenever values like this are quoted there is an implicit
assumption that there is a universal scale on which evidence can be measured.
Currently, there are no developments that support the existence of such a
scale; see the discussion of the Bayes factor in Section 3.4. For the
estimation is natural to quote $\theta_{MLE}(x)$ and report a likelihood
region such as $\{\theta:L^{rel}(\theta\,|\,x)\geq1/r\},$ for some $r,$ as an
assessment of the accuracy of $\theta_{MLE}(x).$

Another serious objection to pure likelihood theory, and to the use of the
likelihood function to determine inferences more generally, arises when
nuisance parameters are present, namely, we want to make inference about
$\psi=\Psi(\theta)$ where $\Psi$ is not 1-1. In general there does not appear
to be a way to define a likelihood function for the parameter of interest
$\psi,$ based on the inference base $I$ only, that is consistent with the
motivation for using likelihood in the first place, namely, $L(\theta\,|\,x)$
is proportional to the probability of the observed data as a function of
$\theta.$ It is common in such contexts to use the \textit{profile likelihood
}
\[
L_{\Psi}^{prof}(\psi\,|\,x)=\sup_{\theta\in\Psi^{-1}\{\psi\}}L(\theta\,|\,x)
\]
as a likelihood function. There are many examples that show that a profile
likelihood is not a likelihood and so such usage is inconsistent with the
basic idea underlying likelihood methods. An alternative to the profile
likelihood is the \textit{integrated likelihood}
\[
L_{\Psi}^{int}(\psi\,|\,x)=\int_{\Psi^{-1}\{\psi\}}L(\theta\,|\,x)\,\Pi
(d\theta\,|\,\psi)
\]
where the $\Pi(\cdot\,|\,\psi)$ are probability measures on the pre-images
$\Psi^{-1}\{\psi\},$ see Berger, Liseo and Wolpert (1999). The integrated
likelihood is a likelihood with respect to the sampling model $\{m(\cdot
\,|\,\psi):\psi\in\Psi(\Theta)\}$ where $m(x\,|\,\psi)=\int_{\Psi^{-1}%
\{\psi\}}f_{\theta}(x)\,\Pi(d\theta\,|\,\psi),$ but this requires adding
$\Pi(\cdot\,|\,\psi)$ to the inference base.\smallskip

\noindent\textbf{Example 4.} \textit{location normal.}

Suppose we want to estimate $\psi=\Psi(\mu)=|\mu|$ so the nuisance parameter
is given by $sgn(\mu).$ The likelihood function is given by $L(\mu
\,|\,x)=\exp(-n(\bar{x}-\mu)^{2}/2\sigma_{0}^{2}).$ Since $(\bar{x}-\mu
)^{2}=(\bar{x}-sgn(\mu)\psi)^{2}$ this is minimized by $sgn(\mu)=1$ when
$\bar{x}>0$ and by $sgn(\mu)=-1$ when $\bar{x}<0.$ Therefore, the profile
likelihood for $\psi$ is $L_{\Psi}^{prof}(\psi\,|\,x)=\exp(-n(|\bar{x}%
|-\psi)^{2}/2\sigma_{0}^{2})$ and this depends on the data only through
$|\bar{x}|.$ To see that $L_{\Psi}^{prof}(\cdot\,|\,x)$ is not a likelihood
function observe that the density of $|\bar{x}|$ is given by%
\begin{equation}
\frac{\sqrt{n}}{\sqrt{2\pi}\sigma_{0}}\left\{  \exp\left(  -\frac{n}%
{2\sigma_{0}^{2}}(|\bar{x}|-sgn(\mu)\psi)^{2}\right)  +\exp\left(  -\frac
{n}{2\sigma_{0}^{2}}(|\bar{x}|+sgn(\mu)\psi)^{2}\right)  \right\}
\label{dens1}%
\end{equation}
which is not proportional to $L_{\Psi}^{prof}(\psi\,|\,x).$ For the integrated
likelihood we need to choose $\Pi(\{1\}\,|\,\psi)=p\in(0,1)$ so
\begin{equation}
L_{\Psi}^{int}(\psi\,|\,x)=p\exp(-n(\bar{x}-\psi)2/2\sigma_{0}^{2}%
)+(1-p)\exp(-n(\bar{x}+\psi)^{2}/2\sigma_{0}^{2}). \label{intlik}%
\end{equation}
Note that (\ref{intlik}) is not obtained by integrating (\ref{dens1}) which is
appropriate because $|\bar{x}|$ is not a minimal sufficient statistic.

The profile and integrated likelihoods are to be used just as the original
likelihood is used for inferences about $\theta$ even though the profile
likelihood is not a likelihood. The profile likelihood leads immediately to
the estimate $|\bar{x}|\ $for $\psi$ while (\ref{intlik}) needs to be
maximized numerically. Although, the functional forms look quite different,
the profile and integrated likelihoods here give almost identical numerical
results. Figure 1 is a plot with $n=2,\bar{x}=1.47$ generated from a $N(2,1)$
distribution and the estimates of $\psi$ are both equal to $1.47.$ There is a
some disagreement in the left tail, so some of the reported intervals will
differ, but as $n$ increases such differences disappear. This agreement is
also quite robust to the choice of $p.$ $\blacksquare$\smallskip%
\begin{figure}
[ptb]
\begin{center}
\includegraphics[
height=2.7518in,
width=2.7596in
]%
{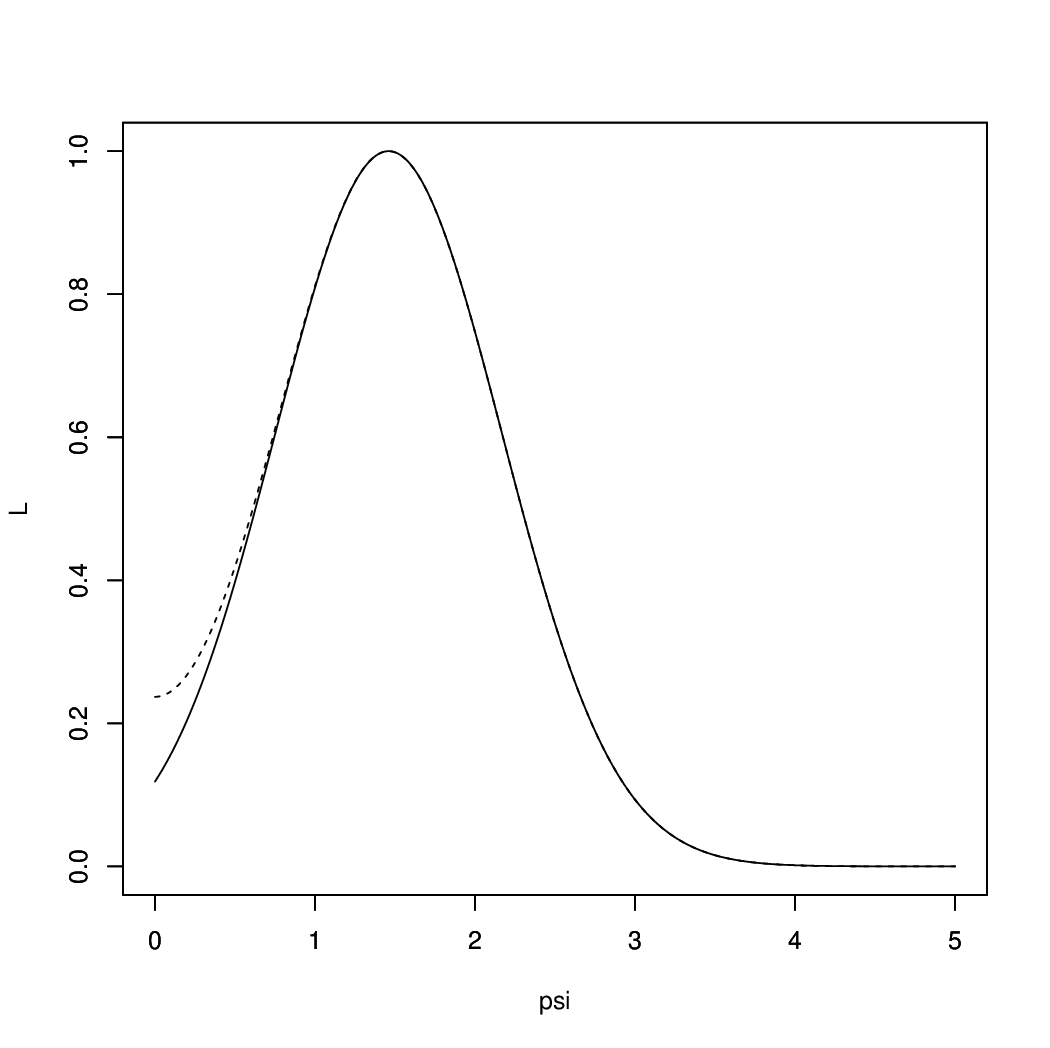}%
\caption{Relative likelihoods for $\psi$ computed from the profile --- and
integrated - - - likelihoods with $n=2,p=1/2$ for data generated from a
$N(2,1)$ distribution in Example 4.}%
\end{center}
\end{figure}

Another issue arises with the profile likelihood, namely, the outcome differs
in contexts which naturally should lead to equivalent results.\smallskip

\noindent\textbf{Example 5.} \textit{prediction with scale normal.}

Suppose $x=(x_{1},\ldots,x_{n})$ is an $iid$ sample from a distribution in the
family $\{n(\cdot,0,\sigma^{2}):\sigma^{2}>0\}$ $.$ Based on the observed
data, the likelihood equals $L(\sigma^{2}\,|\,x)=\sigma^{-n}\exp\{-s_{x}%
^{2}/2\sigma^{2}\}$ where $s_{x}^{2}=\sum_{i=1}^{n}x_{i}{}^{2}$ and the MLE of
$\sigma^{2}$ is $s_{x}^{2}/n.$ Now suppose, however, that interest is in
predicting $k$ future (or occurred but concealed) independent values
$y=(y_{1},\ldots,y_{k})\sim N(0,\sigma^{2}).$ Perhaps the logical predictive
likelihood to use is
\[
L(\sigma^{2},y\,|\,x)=\sigma^{-(n+k)}\exp\{-(s_{x}^{2}+s_{y}^{2})/2\sigma
^{2}\}
\]
where $s_{y}^{2}=\sum_{i=1}^{k}y_{i}^{2}.$ Profiling out $\sigma^{2}$ leads to
the profile MLE of $y$ equaling $0$ for all $x,$ as might be expected.
Profiling $y$ out of $L(\sigma^{2},y\,|\,x),$ however, leads to the profile
likelihood
\[
L^{prof}(\sigma^{2}\,|\,x)=\sigma^{-(n+k)}\exp\{-s_{x}^{2}/2\sigma^{2}\}\neq
L(\sigma^{2}\,|\,x)
\]
for $\sigma^{2}$ and so the profile MLE of $\sigma^{2}$ equals $s^{2}/(n+k).$

When interest is in $\sigma^{2},$ the integrated likelihood for $\sigma^{2}$
produces $L(\sigma^{2}\,|\,x).$ When interest is in $y,\,$then integrating out
$\sigma^{2},$ after placing a probability distribution on $\sigma^{2},$
produces the MLE $0$ for $y$ although the form of the integrated likelihood
will depend on the particular distribution chosen. $\blacksquare$\smallskip

The profile and integrated likelihoods certainly don't always lead to roughly
equivalent results as in Example 4, particularly as the dimension of the
nuisance parameters rises. While the profile likelihood has the advantage of
not having to specify $\Pi(\cdot\,|\,\psi),$ it suffers from a lack of a
complete justification, at least in terms of the likelihood principle and, as
Example 5 demonstrates, it can produce unnatural results. In Section 3.4 it is
shown that the integrated likelihood $m(x\,|\,\psi)$ arises from a very
natural principle characterizing evidence$.$

The discussion here has been mostly about pure likelihood theory where
repeated sampling does not play a role. The likelihood function is also an
important aspect of frequentist inference. Such usage, however, does not lead
to a resolution of the evidence problem, namely, provide a proper
characterization of statistical evidence. In fact, frequentist likelihood
methods use the p-value for \textbf{H}. There is no question that the
likelihood contains the evidence in the data about $\theta,$ but questions
remains as to how to characterize that evidence, whether in favor of or
against a particular value, and also how to express the strength of the
evidence. More discussion of the likelihood principle can be found in Berger
and Wolpert (1988).

\subsection{Bayes Factors and Bayesian Measures of Statistical Evidence}

Bayesian inference adds another ingredient to the inference base, namely, the
\textit{prior} probability measure $\Pi$ on $\Theta,$ so now $I=(\Pi
,\{f_{\theta}:\theta\in\Theta\},x).$ The prior represents beliefs about the
true value of $\theta.$ Note that $(\Pi,\{f_{\theta}:\theta\in\Theta\})$ is
equivalent to a joint distribution for $(\theta,x)$ with density $\pi
(\theta)f_{\theta}(x).$ Once the data $x$ has been observed, a basic principle
of inference is then invoked.

\begin{quote}
Principle of Conditional Probability: for probability model $(\Omega
,\mathcal{F},P),$ if $C\in\mathcal{F}$ is observed to be true, where $P(C)>0,$
then the initial belief that $A\in\mathcal{F}$ is true, as given by $P(A),$ is
replaced by the conditional probability $P(A\,|\,C).$
\end{quote}

\noindent So, we replace the prior $\pi$ by the posterior $\pi(\theta
\,|\,x)=\pi(\theta)f_{\theta}(x)/m(x),$ where $m(x)=\int_{\Theta}\pi
(\theta)f_{\theta}(x)\,d\theta$ is the \textit{prior predictive} density of
the data $x,$ to represent beliefs about $\theta.$ While at times the
posterior is taken to represent the evidence about $\theta,$ this confounds
two distinct concepts, namely, beliefs and evidence. It is clear, however,
that the evidence in the data is what has changed our beliefs and it is this
change that leads to the proper characterization of statistical evidence
through the following principle.

\begin{quote}
Principle of Evidence: for probability model $(\Omega,\mathcal{F},P),$ if
$C\in\mathcal{F}$ is observed to be true where $P(C)>0,$ then there is
evidence in favor of $A\in\mathcal{F}$ being true if $P(A\,|\,C)>P(A),$
evidence against $A\in\mathcal{F}$ being true if $P(A\,|\,C)<P(A)$ and no
evidence either way if $P(A\,|\,C)=P(A).$
\end{quote}

\noindent Therefore, in the Bayesian context we compare $\pi(\theta\,|\,x)$
with $\pi(\theta)$ to determine whether or not there is evidence one way or
the other concerning $\theta$ being the true value. This seems immediate in
the context where $\Pi$ is a discrete probability measure. It is also relevant
in the continuous case, where densities are defined via limits, as in
$\pi(\theta)=\lim_{\epsilon\rightarrow0}\Pi(N_{\epsilon}(\theta))/\nu
(N_{\epsilon}(\theta))$ where $N_{\epsilon}(\theta)$ is a sequence of
neighborhoods of $\theta$\ converging nicely to $\theta$ as $\epsilon
\rightarrow0$ and $\Pi$ is absolutely continuous with respect to support
measure $\nu.$ This leads to the usual expressions for densities, see Evans
(2015), Appendix A for details.

The Bayesian formulation has a very satisfying consistency property. If
interest is in the parameter $\psi=\Psi(\theta),$ then the nuisance parameters
can be integrated out using the conditional prior $\Pi(\cdot\,|\,\psi)$ and
the inference base $I=(\Pi,\{f_{\theta}:\theta\in\Theta\},x)$ is replaced by
$I_{\Psi}=(\Pi_{\Psi},\{m(\cdot\,|\,\psi):\psi\in\Psi(\Theta)\},x)$ where
$\Pi_{\Psi}$ is the marginal prior on $\psi$, with density $\pi_{\Psi}.$
Applying the principal of evidence here means to compare the posterior%
\begin{equation}
\pi_{\Psi}(\psi\,|\,x)=\pi_{\Psi}(\psi)\frac{m(x\,|\,\psi)}{m(x)}=\int
_{\Psi^{-1}\{\psi\}}\Pi(d\theta\,|\,x) \label{post}%
\end{equation}
to the prior density $\pi_{\Psi}(\psi)=\int_{\Psi^{-1}\{\psi\}}\Pi(d\theta).$
So, there is evidence in favor of $\psi$ being the true value whenever
$\pi_{\Psi}(\psi\,|\,x)>\pi_{\Psi}(\psi),$ etc.

In many applications we need more than the simple characterization of evidence
that the principle of evidence gives us. A \textit{valid measure of evidence}
is then any real-valued function of $I_{\Psi}$ that satisfies the existence of
a cut-off $c$ such that the function taking a value greater than $c$
corresponds to evidence in favor, etc. One very simple function satisfying
this is the \textit{relative belief ratio}
\begin{equation}
RB_{\Psi}(\psi\,|\,x)=\frac{\pi_{\Psi}(\psi\,|\,x)}{\pi_{\Psi}(\psi)}%
=\frac{m(x\,|\,\psi)}{m(x)} \label{rb1}%
\end{equation}
where the second equality follows from (\ref{post}) and is sometimes referred
to as the Savage-Dickey ratio, see Dickey (1971). Using $RB_{\Psi}%
(\psi\,|\,x)$ the values of $\psi$ are now totally ordered with respect to
evidence, as when $1<RB_{\Psi}(\psi_{1}\,|\,x)<RB_{\Psi}(\psi_{2}\,|\,x)$
there is more evidence in favor of $\psi_{2}$ than for $\psi_{1},$ etc.

This suggests an immediate answer to \textbf{E}, namely, record the
\textit{relative belief estimate} $\psi(x)=\arg\sup_{\psi}RB_{\Psi}%
(\psi\,|\,x)$ as this value has the maximum evidence in its favor. Also, to
assess the accuracy of $\psi(x)$ record the \textit{plausible region}
$Pl_{\Psi}(x)=\{\psi:RB_{\Psi}(\psi\,|\,x)>1\},$ the set of $\psi$ values with
evidence in their favor, together with its posterior content $\Pi_{\Psi
}(Pl_{\Psi}(x)\,|\,x)$ as this measures the belief that the true value is in
$Pl_{\Psi}(x).$ Note that while $\psi(x)$ depends on the choice of the
relative belief ratio to measure evidence, the plausible region $Pl_{\Psi}(x)$
only depends on on the principle of evidence. This suggests that any other
valid measure of evidence can be used instead to determine an estimate as it
will lie in $Pl_{\Psi}(x)$. A $\gamma$-credible region $C_{\Psi,\gamma}(x)$
for $\psi$ can also be quoted based on any valid measure of evidence provided
$C_{\Psi,\gamma}(x)\subset Pl_{\Psi}(x)$ as otherwise $C_{\Psi,\gamma}(x)$
would contain a value $\psi$ for which there is evidence against $\psi$ being
the true value. For example a relative belief $\gamma$-credible region takes
the form $C_{\Psi,\gamma}(x)=\{\psi:RB_{\Psi}(\psi\,|\,x)\geq c_{\gamma}(x)\}$
where $c_{\gamma}(x)=\inf\{c:\Pi_{\Psi}(RB_{\Psi}(\psi\,|\,x)\leq
c\,|\,x)\geq1-\gamma\}.$

For \textbf{H}, the evidence concerning $H_{0}:\Psi(\theta)=\psi_{0},$ then
the value $RB_{\Psi}(\psi_{0}\,|\,x)$ tells us immediately whether there is
evidence in favor of or against $H_{0}.$ To measure the strength of the
evidence concerning $H_{0}$ there is the posterior probability $\Pi_{\Psi
}(\{\psi_{0}\}\,|\,x),$ as this measures our belief in what the evidence says.
For if $RB_{\Psi}(\psi_{0}\,|\,x)>1$ and $\Pi_{\Psi}(\{\psi_{0}%
\}\,|\,x)\approx1,$ then there is strong evidence that $H_{0}$ is true while
when $RB_{\Psi}(\psi_{0}\,|\,x)<1$ and $\Pi_{\Psi}(\{\psi_{0}\}\,|\,x)\approx
0,$ then there is strong evidence that $H_{0}$ is false. Often, however,
$\Pi_{\Psi}(\{\psi_{0}\}\,|\,x)$ will be small, even 0 in the continuous case,
so it makes more sense to measure the strength of the evidence in such a case
by $Str_{\Psi}(\psi_{0})=\Pi_{\Psi}(RB_{\Psi}(\psi\,|\,x)\leq RB_{\Psi}%
(\psi_{0}\,|\,x)\,|\,x)$ as when $RB_{\Psi}(\psi_{0}\,|\,x)>1$ and $Str_{\Psi
}(\psi_{0})\approx1,$ then there is small belief that the true value of $\psi$
has more evidence in its favor than $\psi_{0},$ etc.

Recalling the discussion in Section 3.1 about the \textit{difference that
matters} $\delta,$ it is relatively easy to take this into account in this
context, at least when $\psi$ is real valued. For this we consider a grid of
values $\ldots\psi_{-2},\psi_{-1},\psi_{0},\psi_{1},\psi_{2},\ldots$ separated
by $\delta$ and then conduct the inferences using the relative belief ratios
of the intervals $[\psi_{i}-\delta/2,\psi_{i}+\delta/2).$ In effect, $H_{0}$
is now $H_{0}:\Psi(\theta)\in\lbrack\psi_{0}-\delta/2,\psi_{0}+\delta
/2).$\smallskip

\noindent\textbf{Example 6.} \textit{location normal.}

Suppose that we add the prior $\mu\sim N(\mu_{0},\tau_{0}^{2})$ to form a
Bayesian inference base. The posterior distribution is then
\[
\mu\,|\,x\sim N\left(  \left(  n/\sigma_{0}^{2}+1/\tau_{0}^{2}\right)
^{-1}(n\bar{x}/\sigma_{0}^{2}+\mu_{0}/\tau_{0}^{2}),\left(  n/\sigma_{0}%
^{2}+1/\tau_{0}^{2}\right)  ^{-1}\right)  .
\]
Suppose our interest is inference about $\psi=\Psi(\mu)=|\mu|.$ Figure 2 is a
plot of the relative belief ratio based on the data in Example 4 $(n=2,\bar
{x}=1.47$ generated from a $N(2,1)$ distribution) and using the prior with
$\mu_{0}=0,\tau_{0}=2.$ From (\ref{rb1}) it is seen that maximizing $RB_{\Psi
}(\psi\,|\,x)$ is equivalent to maximizing the integrated likelihood and in
this case the prior is such that $sgn(\mu)=1$ with prior probability 0.5 so
indeed the relative belief estimate is $\psi(x)=1.47.$ To assess the accuracy
of the estimate, we record the plausible interval $Pl_{\Psi}(x)=(0.65,2.26)$
which has posterior content $0.76.$ There is evidence in favor of $H_{0}%
:\psi=2,$ since $2\in Pl_{\Psi}(x),$ and $RB_{\Psi}(2\,|\,x)=1.41$ with
$Str_{\Psi}(\psi_{0})=0.42,$ so there is only moderate evidence in favor of
$2$ being the true value of $\psi$. Of course, these results improve with
sample size $n,$ For example, for $n=10,\psi(x)=1.83$ with $Pl_{\Psi
}(x)=(1.40,2.27)$ with posterior content $0.95$ and $RB_{\Psi}(2\,|\,x)=5.00$
with strength $0.44.$ So, the plausible interval has shortened considerably
and its posterior content increased but the strength of the evidence in favor
of $H_{0}:\psi=2$ has not increased by much. It is important to note that
because we are employing a discretization (with $\delta=0.01$ here), and since
the posterior is inevitably virtually completely concentrated in $[\psi
_{0}-\delta/2,\psi_{0}+\delta/2)$ as $n\rightarrow\infty,$ the strength
$Str_{\Psi}(2)$ converges to 1 at any the values in this interval and to 0 at
values outside. $\blacksquare$
\begin{figure}
[ptb]
\begin{center}
\includegraphics[
height=2.8193in,
width=2.8271in
]%
{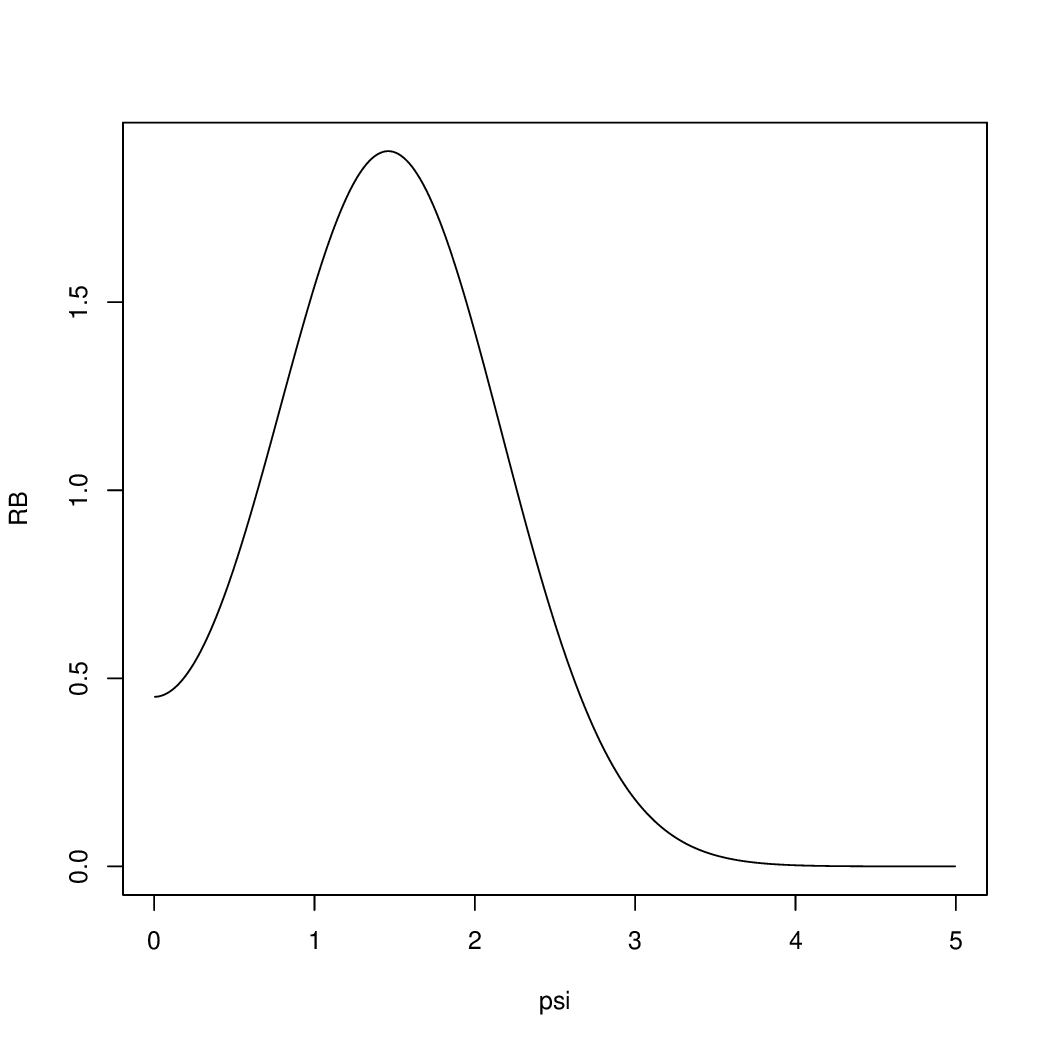}%
\caption{Relative belief ratio for $\psi$ in Example 5 with $n=2$.}%
\end{center}
\end{figure}

There are a number of optimal properties for the relative belief inferences as
discussed in Evans (2015). While this is perhaps the first attempt to build a
theory of inference based on the principle of evidence, there is a
considerable literature on this idea in the philosophy of science where it
underlies what is known as confirmation theory, see Salmon (1973). For
example, Popper (1968) (Appendix ix) writes

\begin{quote}
If we are asked to give a criterion of the fact that the evidence $y$ supports
or corroborates a statement $x$, the most obvious reply is: that $y$ increases
the probability of $x.$
\end{quote}

\noindent One issue with the philosophers' discussions is that these are never
cast in a statistical context. Many of the anomalies raised in those
discussions, such as Hempel's the Raven Paradox, can be resolved when
formulated as statistical problems. Also, the relative belief ratio itself has
appeared elsewhere, although under different names, as a natural measure of evidence.

There are indeed other valid measures of evidence besides the relative belief
ratio. For example, the Bayes factor originated with Jeffreys (1935), see Kass
and Raftery (1995) for discussion. For probability model $(\Omega
,\mathcal{F},P),A,C\in\mathcal{F}$ with $P(A)>0,P(C)>0,$ the \textit{Bayes
factor} in favor of $A$ being true, after observing that $C$ is true, is given
by
\[
BF(A\,|\,C)=\frac{P(A\,|\,C)/P(A^{c}\,|\,C)}{P(A)/P(A^{c})},
\]
the ratio of the posterior odds in favor of $A\,$to the prior odds in favor of
$A.$ It is immediate that $BF(A\,|\,C)$ is a valid measure of evidence because
$BF(A\,|\,C)>1$ iff $P(A\,|\,C)>P(A).$ Also, $BF(A\,|\,C)=RB(A\,|\,C)/RB(A^{c}%
\,|\,C)$\ and so, because $RB(A\,|\,C)>1$ iff $RB(A^{c}\,|\,C)<1,$ this
implies $BF(A\,|\,C)>RB(A\,|\,C)$ when there is evidence in favor and
$BF(A\,|\,C)<RB(A\,|\,C)$ when there is evidence against. These inequalities
are important because it is sometimes asserted that the Bayes factor is not
only a measure of evidence but its value is a measure of the strength of that
evidence. Table \ref{tab1} gives a scale, due to Jeffreys (1961) (Appendix B),
which supposedly can be used to assess the strength of evidence as given by
the Bayes factor. Again, this is an attempt to establish a universal scale on
which evidence can be measured and currently no grounds exist for claiming
that such a scale exists. If we consider that the strength of the evidence is
to be measured by how strongly we believe what the evidence says, then simple
examples can be constructed to show that such a scale is inappropriate.
Effectively, the strength of the evidence is context dependent and needs to be
calibrated.\smallskip

\noindent\textbf{Example 7.}

Suppose $\Omega$ is a set with $\#(\Omega)=N.$ A value $\omega$ is generated
uniformly from $\Omega$ and partially concealed but it is noted that
$\omega\in C$ where $\#(C)=n.$ It is desired to know if $\omega\in A$ where
$\#(A)=1$ and $A\subset C.$ Then we have that $P(A)=1/N$ and $P(A\,|\,C)=1/n$
so $RB(A\,|\,C)=N/n>1$ since\thinspace$n<N.\,$\ The posterior belief in what
the evidence is saying is, however, $P(A\,|\,C)=1/n$ which can be very small,
say $n=10^{3}.$ But, if $N=10^{6},$ then $RB(A\,|\,C)=10^{3}$ which is well
into the range where Jeffreys scale says it is decisive evidence in favor of
$A.$ Clearly there is only very weak evidence in favor of $A$ being true here.
Also, note that the posterior probability by itself does not indicate evidence
in favor although the observation that $C$ is true must be evidence in favor
of $A$ because $A\subset C.$ The example can obviously be modified to show
that any such scale is not appropriate. $\blacksquare$

.
\begin{table}[tbp] \centering
\begin{tabular}
[c]{|l|l|}\hline
$BF$ & Strength\\\hline
$1$ to $10^{1/2}$ & Barely worth mentioning\\\hline
$10^{1/2}$ to $10$ & Substantial\\\hline
$10$ to $10^{3/2}$ & Strong\\\hline
$10^{3/2}$ to $10^{2}$ & Very Strong\\\hline
$>$ $10^{2}$ & Decisive\\\hline
\end{tabular}
\caption{Jeffreys' Bayes factor scale for measuring strength of the evidence in favor (the strength of evidence against is measured by the reciprocals),}\label{tab1}%
\end{table}%

There is another issue with current usage of the Bayes factor that needs to be
addressed. This arises when interest is in assessing the evidence for or
against $H_{0}:\Psi(\theta)=\psi_{0}$ when $\Pi_{\Psi}(\{\psi_{0}\})=0.$
Clearly the Bayes factor is not defined in such a context. An apparent
solution is provided by choosing a prior of the form $\Pi_{H_{0},p}%
=p\Pi_{H_{0}}+(1-p)\Pi$ where $p>0$ is a prior probability assigned to $H_{0}$
and $\Pi_{H_{0}}$ is a prior probability measure on the set $\Psi^{-1}%
\{\psi_{0}\}.$ Since $\Pi_{p}(\Psi^{-1}\{\psi_{0}\})=p$ the Bayes factor is
now defined. Of some concern now is how $\Pi_{H_{0}}$ should be chosen. A
natural choice is to take $\Pi_{H_{0}}=\Pi(\cdot\,|\,\psi_{0})$ otherwise
there is a contradiction between beliefs as expressed by $\Pi$ and $\Pi
_{H_{0}}.$ It follows very simply, however, that in the continuous case, when
$\Pi_{H_{0}}=\Pi(\cdot\,|\,\psi_{0}),$ then $BF(\psi_{0}\,|\,x)=RB_{\Psi}%
(\psi_{0}\,|\,x).$ Moreover, if instead of using such a mixture prior we
instead define the Bayes factor via
\[
BF_{\Psi}(\psi_{0}\,|\,x)=\lim_{\epsilon\rightarrow0}BF(N_{\epsilon}(\psi
_{0})\,|\,x)\lim_{\epsilon\rightarrow0}\frac{\Pi_{\Psi}(N_{\epsilon}(\psi
_{0})\,|\,x)/\Pi_{\Psi}(N_{\epsilon}^{c}(\psi_{0})\,|\,x)}{\Pi_{\Psi
}(N_{\epsilon}(\psi_{0}))/\Pi_{\Psi}(N_{\epsilon}^{c}(\psi_{0}))},
\]
where $N_{\epsilon}(\psi_{0})$ is a sequence of neighborhoods of $\psi_{0}$
converging nicely to $\psi_{0},$ then, under weak conditions ($\pi_{\Psi}$ is
positive and continuous at $\psi_{0}),$ $BF_{\Psi}(\psi_{0}\,|\,x)=RB_{\Psi
}(\psi_{0}\,|\,x).$ So, there is no need to introduce the mixture prior to get
a valid measure of evidence. Furthermore, now the Bayes factor is available
for \textbf{E }as well as \textbf{H}, as surely any proper characterization of
statistical evidence must be, and the relevant posterior for $\psi$ is
$\Pi_{\Psi}(\cdot\,|\,x)$ as obtained from $\Pi$ rather than from $\Pi
_{H_{0},p}.$ There are many other reasons to prefer the relative belief ratio
to the Bayes factor as a measure of evidence and these are discussed in
Al-Labadi, Alzaatreh and Evans (2024).

Other Bayesian measures of evidence have been proposed. For example, Pereira
and Stern (1999) proposed to measure the evidence for $H_{0}:\Psi(\theta
)=\psi_{0}$ by computing $\theta_{\ast}=\arg\sup_{\Psi(\theta)=\psi_{0}}%
\pi(\theta\,|\,x)$ and then use the posterior tail probability $ev(H_{0}%
\,|\,x)=\Pi(\pi(\theta\,|\,x)\leq\pi(\theta_{\ast}\,|\,x)\,|\,x).$ If
$ev(H_{0}\,|\,x)$ is large, this is evidence against $H_{0}$ while if it is
small, it is evidence in favor. Note that $ev(H_{0}\,|\,x)$ is sometimes
referred to as an e-value but this is different than the e-values discussed in
Section 3.1. Further discussion and development of this concept can be found
in\ Stern and Pereira (2022). Clearly this is building on the idea that
underlies the p-value, namely, providing a measure that locates a point in a
distribution and using this to assess evidence. It does not, however, conform
to the principle of evidence.

\subsubsection{Bias and Error Probabilities}

One criticism that is made of Bayesian inference is that there are no measures
of reliability of the inferences as is an inherent part of frequentism. It is
natural to add such measures, however, to assess whether or not the specified
model and prior could potentially lead to misleading inferences. For example,
suppose it could be shown that evidence in favor of $H_{0}:\Psi(\theta
)=\psi_{0}$ would be obtained with prior probability near 1, for a data set of
given size. It seems obvious that, if we did collect this amount of data and
obtained evidence in favor of $H_{0},$ then this fact would undermine our
confidence in the reliability of the reported inference. \smallskip

\noindent\textbf{Example 8.} \textit{location nornal and the Jeffreys-Lindley
paradox.}

Suppose we have the location normal model $\bar{x}\sim N(\mu,\sigma_{0}%
^{2}/n)$ and $T_{\mu_{0}}(x)=\sqrt{n}|\bar{x}-\mu_{0}|/\sigma_{0}=5$ is
obtained which leads to the p-value $2(1-\Phi(\sqrt{n}|\bar{x}-\mu_{0}%
|/\sigma_{0}))=5.733031\times10^{-07}$ so there would appear to be definitive
evidence that $H_{0}:\mu=\mu_{0}$ is false. Suppose the prior $\mu\sim
N(\mu_{0},\tau_{0}^{2})$ is used and the analyst chooses $\tau_{0}^{2}$ very
large to reflect the fact that little is known about the true value of $\mu.$
It can be shown (see Evans (2015)), however, that $RB_{\Psi}(\mu
_{0}\,|\,x)\rightarrow\infty$ as $\tau_{0}^{2}\rightarrow\infty.$ Therefore,
for a very diffuse prior, evidence in favor of $H_{0}$ will be obtained and so
the frequentist and Bayesian will disagree. Note that the Bayes factor equals
the relative belief ratio in this context. A partial resolution of this
contradiction is obtained by noting that $Str_{\Psi}(\mu_{0})\rightarrow
2(1-\Phi(\sqrt{n}|\bar{x}-\mu_{0}|/\sigma_{0}))$ as $\tau_{0}^{2}%
\rightarrow\infty$ and so the Bayesian measure of evidence is only providing
very weak evidence when the p-value is small. This anomaly occurs even when
the true value of $\mu$ is indeed far away from $\mu_{0}$ and so the fault
here does not lie with the p-value. $\blacksquare$\smallskip

Prior error probabilities associated with Bayesian measures of evidence can,
however, be computed and these lead to a general resolution of the apparent
paradox. There are two error probabilities for \textbf{H}\ that we refer to as
\textit{bias against} $H_{0}$ and \textit{bias in favor of} $H_{0},$ as given
by the following two prior probabilities%

\begin{align*}
&  \text{\textit{bias against}}_{\Psi}(\psi_{0})=M(RB_{\Psi}(\psi
_{0}\,|\,X)\leq1\,|\,\psi_{0}),\\
&  \text{\textit{bias in favor of}}_{\Psi}(\psi_{0})=\sup_{\psi:d_{\Psi}%
(\psi,\psi_{0})\geq\delta}M(RB_{\Psi}(\psi_{0}\,|\,X)\geq1\,|\,\psi).
\end{align*}
Both of these biases are independent of the valid measure of evidence used as
they only depend on the principle of evidence as applied to the model and
prior chosen. The \textit{bias against}$_{\Psi}(\psi_{0})$ is the prior
probability of not getting evidence in favor of $\psi_{0}$ when it is true and
plays a role similar to type I error. The \textit{bias in favor of}$_{\Psi
}(\psi_{0})$ is the prior probability of not getting evidence against
$\psi_{0}$ when it is meaningfully false and plays a role similar to type II
error. Here $\delta$ is the deviation from $\psi_{0}$ which is of scientific
interest as determined by some measure of distance $d_{\Psi}$ on $\Psi
(\Theta).$ As discussed in Evans (2015), these biases cannot be controlled by,
for example, the choice of prior as a prior that reduces \textit{bias
against}$(\psi_{0})$ causes bias in favor of $(\psi_{0})$ to increase and
conversely. The proper control of these quantities is through sample size as
both converge to 0 as $n\rightarrow\infty.$\smallskip

\noindent\textbf{Example 9.} \textit{location nornal and the Jeffreys-Lindley
paradox.}

In this case \textit{bias against}$(\mu_{0})\rightarrow0$ and the \textit{bias
in favor of }$(\mu_{0})\rightarrow1$ as $\tau_{0}^{2}\rightarrow\infty.$ This
leads to an apparent resolution of the issue: don't choose the prior to be
arbitrarily diffuse to reflect noninformativeness, rather choose a prior that
is sufficiently diffuse to cover the interval where it is known $\mu$ must
lie, e.g., the interval where the measurements lie, and then choose $n$ so
that both biases are suitably small. While the Jeffreys-Lindley paradox arises
due to diffuse priors inducing bias in favor, an overly concentrated prior
induces bias against, but again this bias can be controlled via the amount of
data collected. $\blacksquare$\smallskip

There are also two biases associated with \textbf{E} obtained by averaging,
using the prior $\Pi_{\Psi},$ the biases for \textbf{H}. These can also be
expressed in terms of coverage probabilities of the plausible region
$Pl_{\Psi}(x).$ The biases for \textbf{E} are given by%
\begin{align*}
&  \text{\textit{bias against}}_{\Psi}=E_{\Pi_{\Psi}}(M(RB_{\Psi}%
(\psi\,|\,X)\leq1\,|\,\psi))=E_{\Pi_{\Psi}}\left(  M(\psi\notin Pl_{\Psi
}(X)\,|\,\psi)\right)  ,\\
&  \text{\textit{bias in favor of}}_{\Psi}=E_{\Pi_{\Psi}}(\sup_{\psi^{\prime
}:d_{\Psi}(\psi^{\prime},\psi)\geq\delta}M(RB_{\Psi}(\psi\,|\,X)\geq
1\,|\,\psi^{\prime}))\\
&  \hspace{1in}=E_{\Pi_{\Psi}}\left(  \sup_{\psi^{\prime}:d_{\Psi}%
(\psi^{\prime},\psi)\geq\delta}M(\psi\notin Im_{\Psi}(X)\,|\,\psi^{\prime
})\right)  ,
\end{align*}
where $Im_{\Psi}(x)=\{\psi:RB_{\Psi}(\psi\,|\,x)<1\}$ is the
\textit{implausible region}, namely, the set of values for which evidence
against has been obtained. So, \textit{bias against}$_{\Psi}$ is the prior
probability that true value is not in the plausible region $Pl_{\Psi}(x)$ and
so $1-E_{\Pi_{\Psi}}\left(  M(\psi\notin Pl_{\Psi}(X)\,|\,\psi)\right)  $ is
the prior coverage probability (\textit{Bayesian confidence}) of $Pl_{\Psi
}(x).$ It is of interest that there is typically a value $\psi_{\ast}=\arg
\sup_{\psi}$ \textit{bias against}$_{\Psi}(\psi)$ and so
\[
M(\psi\in Pl_{\Psi}(X)\,|\,\psi)\geq1-M(RB_{\Psi}(\psi_{\ast}\,|\,X)\leq
1\,|\,\psi_{\ast})
\]
gives a lower bound on the frequentist confidence of $Pl_{\Psi}(x)$ with
respect to the model $\{m(\cdot\,|\,\psi):\psi\in\Psi\}$\ obtained from the
original model by integrating out the nuisance parameters. In both cases the
Bayesian confidences are average frequentist confidences with respect to the
original model. The \textit{bias in favor of}$_{\Psi}$ is the prior
probability that a meaningfully false value is not in the implausible region.
Again, both of these biases do not depend on the valid measure of statistical
evidence used and converge to 0 with increasing sample size.

With the addition of the biases, a link is established between Bayesianism and
frequentism: inferences are Bayesian but the reliability of the inferences is
assessed via frequentist criteria. More discussion on bias can be found in
Evans and Guo (2021) and Evans et al. (2024) demonstrates how the biases can
be used to resolve some issues that arise with the confidence concept in general.

\section{Conclusions}

This paper has been a review of various conceptions of statistical evidence as
this has been discussed in the literature. It is seen that precisely
characterizing evidence is still an unsolved problem in purely frequentist
inference and in pure likelihood theory but the Bayesian framework provides a
natural way to do this through the principle of evidence. One often hears
complaints about the need to choose a prior which can be done properly through
elicitation, see O'Hagan et al. (2006). Furthermore, priors are falsifiable as
in checking for prior-data conflict, see Evans and Moshonov (2006) and Nott et
al. (2020). Finally, the effects of priors can be mitigated by the control of
the biases which are directly related to the principle of evidence.

\section{References}

\noindent Al-Labadi, L., Alzaatreh, A. and Evans, M. (2024) How to measure
evidence and its strength: Bayes factors or relative belief ratios?
arXiv:2301.08994\smallskip

\noindent Berger, J.O. and Wolpert, R.L. (1988) The likelihood principle: A
review, generalizations, and statistical implications. IMS Lecture Notes
Monogr. Ser., 6.\smallskip

\noindent Berger, J.O., Liseo, B. and Wolpert, R.L. (1999) Integrated
Likelihood Methods for Eliminating Nuisance Parameters. Statistical Science,
14, 1, 1--22.\smallskip

\noindent Birnbaum, A. (1962) On the foundations of statistical inference
(with discussion). Journal of the American Statistical Association. 57 (298),
269--326.\smallskip

\noindent Birnbaum, A. (1977) The Neyman-Pearson theory as decision theory,
and as inference Theory; with criticism of the Lindley-SavageaArgument for
Bayesian theory. Synthese 36, 19-49.\smallskip

\noindent Boring, E.G. (1919) Mathematical vs. Scientific Significance.
Psychological Bulletin, 16, 335-338.\smallskip

\noindent Cornfield, J. Sequential Trials, Sequential Analysis and the
Likelihood Principle. The American Statistician, 20, 2, 18--23.\smallskip

\noindent Dickey, J. (1971). The weighted likelihood ratio, linear hypotheses
on normal location parameters. Annals of Statistics, 42, 204--223.\smallskip

\noindent Edwards, A.W.F. (1992). Likelihood (2nd ed.). Baltimore, MD: Johns
Hopkins University Press.\smallskip

\noindent Evans, M., Fraser, D.A.S. and Monette, G. (1986) On principles and
arguments to likelihood - with discussion. Canad. J. of Statist., 14, 3,
181-199.\smallskip

\noindent Evans, M. and Moshonov, H. (2006) Checking for prior-data conflict.
Bayesian Analysis, 1, 4, 893-914.\smallskip\ 

\noindent Evans, M. (2013) What does the proof of Birnbaum's theorem prove?
Electronic Journal of Statistics, 7, 2645-2655.\smallskip

\noindent Evans, M. (2015) Measuring Statistical Evidence Using Relative
Belief. Chapman \& Hall/CRC Monographs on Statistics \& Applied
Probability.\smallskip

\noindent Evans, M. and Guo, Y. (2021) Measuring and controlling bias for some
Bayesian inferences and the relation to frequentist criteria. Entropy, 23(2),
190, doi: 10.3390/e23020190.\smallskip

\noindent Evans, M. Liu, M. Moon, M., Sixta, S. Wei, S. and Yang, S. (2024) On
some problems of Bayesian region construction with guaranteed coverages.
Statistical Papers, 65:309--334.\smallskip

\noindent Evans, M., Frangakis, C. (2023) On resolving problems with
conditionality and its implications for characterizing statistical evidence.
Sankhya A 85, 1103--1126.\smallskip

\noindent Fisher R. A. (1922) On the mathematical foundations of theoretical
statistics. Philosophical Transactions of the Royal Society of London. Series
A, 222, 309--368.

\noindent Fisher, R.A. (1925) Statistical Methods for Research Workers (14th
edition). Hafner Press.\smallskip\ 

\noindent Giere, R. N. (1977) Publications by Allan Birnbaum. Synthese, 36, 1,
15--17.\smallskip

\noindent Grunwald, P., de Heide, R. and Koolen, W. M. (2024). Safe testing.
To appear in Journal of the Royal Statistical Society, Series B.\smallskip

\noindent Jeffreys H. (1935) Some tests of significance, treated by the theory
of probability. Mathematical Proceedings of the Cambridge Philosophical
Society, 31(2):203-222.\smallskip

\noindent Jeffreys H. (1961) Theory of Probability (Third Edition).
Oxford.\smallskip

\noindent Kass, R. E. and Raftery, A. E. (1995). Bayes Factors. Journal of the
American Statistical Association, 90(430), 773--795.\smallskip

\noindent Nott,D., Wang, X., Evans, M., and Englert, B-G. (2020) Checking for
prior-data conflict using prior to posterior divergences. Statistical Science,
35, 2, 234-253.\smallskip

\noindent O'Hagan, A., Buck, C. E., Daneshkhah, A., Eiser, J. R., Garthwaite,
P. H., Jenkinson, D. J., Oakley, J. E., and Rakow, T. (2006), Uncertain
Judgements: Eliciting Expert Probabilities. Wiley.\smallskip

\noindent Pereira, C. A. de B. and Stern, J. M. (1999) Evidence and
credibility: full Bayesian significance test for precise hypotheses. Entropy,
1, 99-110.\smallskip

\noindent Popper, K. (1968) The Logic of Scientific Discovery, Harper
Torchbooks.\smallskip

\noindent Royall, Richard M. (1997). Statistical Evidence: A likelihood
paradigm. London, UK: Chapman \& Hall.\smallskip

\noindent Salmon,W. (1973) Confirmation. Scientific American, 228,
75--81.\smallskip

\noindent Shafer, G., Shen, A., Vereshchagin, N. and Vovk, V. (2011) Test
martingales, Bayes factors and p-values. Statistical Science, 26, (1), 84 -
101.\smallskip

\noindent Stern, J. M. and Pereira, C. A. de B. (2022) The e-value and the
full Bayesian significance test: logical properties and philosophical
consequences. S\~{a}o Paulo Journal of Mathematical Sciences,16 (1),
566-584.\smallskip

\noindent Stigler, S.M. (1986)\ The History of Statistics. Belknap
Press.\smallskip

\noindent Vovk, V. and Wang, R. (2023). Confidence and discoveries with
e-values. Statistical Science, 38(2), 329-354.
\end{document}